\documentclass[a4paper,12pt]{exam}
\printanswers 
\qformat{\textbf{Problema \thequestion}\hfill}
\usepackage{color} 
\definecolor{SolutionColor}{rgb}{0.8,0.9,1} 
\shadedsolutions 

\footer{}{Pag.~\thepage}{}

\usepackage{graphicx}
\usepackage{amsmath, graphicx}

\usepackage{pgf}
\usepackage{tikz,wrapfig}
\usetikzlibrary{arrows}

\usepackage[utf8]{inputenc}
\usepackage{amsfonts, amsthm, amsmath, enumerate}
\usepackage{mathrsfs}
\usepackage[pdftex,pdfauthor={Dario Trevisan}]{hyperref}
\hypersetup{colorlinks=true}

\newtheorem{theorem}{Theorem}
\newtheorem{proposition}[theorem]{Proposition}
\newtheorem{lemma}[theorem]{Lemma}
\newtheorem{corollary}[theorem]{Corollary}
\theoremstyle{definition}

\theoremstyle{remark}

\newtheorem{remark}{Remark}

\usepackage{amsmath,color,ulem}
\newcommand{\stkout}[1]{\ifmmode\text{\sout{\ensuremath{#1}}}\else\sout{#1}\fi}

\begin{document}

	\title{ On the limit distribution of  extremes of generalized Oppenheim random variables  }
	\author{Rita Giuliano \footnote{Dipartimento di
			Matematica, Universit\`a di Pisa, Largo Bruno
			Pontecorvo 5, I-56127 Pisa, Italy (email: rita.giuliano@unipi.it)}~~and Milto Hadjikyriakou\footnote{School of Sciences, University of Central Lancashire, Cyprus campus, 12-14 University Avenue, Pyla, 7080 Larnaka, Cyprus (email:
			mhadjikyriakou@uclan.ac.uk).  }} 	
	
	\maketitle
	
	\begin{abstract}
		This paper investigates the asymptotic behavior of the  extremes of a sequence   of   generalized Oppenheim random variables. Particularly, we establish conditions under which some  normalized extremes of sequences  arising from Oppenheim expansions belong to the maximum domain of attraction of the Fr\'echet distribution. Additionally, we identify conditions under which the maxima and minima of Oppenheim random variables demonstrate some kind of asymptotic independence. Finally, we prove an Extreme Types theorem for Oppenheim expansions with unknown dependent structure.
	\end{abstract}
	
	
	\medskip
	
	

	
	\section{Introduction} 
	One of the main topics in the so-called \textit{Extreme Value Theory} is the study of the asymptotic behaviour of
	$$M_n = \max\{U_1,\ldots,U_n\}$$
	of $n$ independent and identically distributed random variables. The celebrated \textit{Extremal Types Theorem} (known also as {\it   Fisher-Tippett-Gnedenko Theorem)}  states that if   there are constants $a_n>0, b_n$ such that
	$$P(a_n(M_n-b_n)\leq x)\to G(x),$$ 
	where $G$ is a non degenerate distribution, then $G$ must be one of the following types:
	\begin{enumerate}
		\item[(I)] Gumbel distribution:$$G(x) = \exp(-e^{-x}),\, \, -\infty<x<\infty;$$ 
		\item[(II)]  Fr\'echet distribution:$$G(x) = \begin{cases}
			\exp(-x^{-\alpha}), &  x>0, \quad (\alpha>0)\\
			0, & x\leq 0;
		\end{cases}$$
		\item[(III)] negative Weibull distribution: $$G(x) = \begin{cases}
			\exp(-(-x)^\alpha), & x\leq 0, \quad (\alpha>0)\\
			1, & x>0. 
		\end{cases} $$   
	\end{enumerate}
	These are known as the three \textit{Extremal Value Distributions}. 
	
	\bigskip
	\noindent
	As it is customary, we say that a sequence of i.i.d. random variables  belongs to the maximum domain of attraction (MDA) of the distribution $G$ if the distribution of the  sequence of the normalized maxima  converges to $G$.  Thus the Extremal Types Theorem says that the given sequence can be in the MDA   of the Gumbel or Fr\'echet or Weibull distribution only.

	\bigskip
	\noindent
	The Extremal Types Theorem was discussed for the first time in 1928 by Fisher and Tippett \cite{FT1928} and its most general form was proven by Gnedenko in 1943 \cite{G1943} while a simpler approach was provided by de Haan in 1976 \cite{H1976}. In 1974, Leadbetter \cite{L1974} showed that the Extremal Types Theorem remains true for random variables that are stationary and satisfy a quite weak dependence structure which is referred as\textit{ Leadbetter's condition} or as the \textit{distributional mixing condition} (see the Discussion Section).
	
	\bigskip
	\noindent
	In this paper we are interested in studying the limit behaviour of the extremes of a special class of random variables known as generalized Oppenheim random variables denoted by $R_n$ and described as follows: let $(B_n)_{n\geq  1}$ be a sequence of integer valued random  variables defined on $(\Omega, \mathcal{A}, P)$, where $\Omega =[0,1]$, $\mathcal{A}$ is the $\sigma$-algebra of the Borel subsets of $[0,1]$ and $P$ is the Lebesgue measure on $[0,1]$. Let $\{F_n, n\geq 1\}$ be a sequence of probability distribution functions with  $F_n(0)=0$, $\forall n$ and moreover let $\varphi_n:\mathbb{N}^*\to \mathbb{R}^+$ be a sequence of functions. Furthermore, let $(q_n)_{n\geq  1}$ with $q_n=q_n(h_1, \dots, h_n)$ be a sequence of nonnegative numbers (i.e. possibly depending on the $n$ integers $h_1, \dots, h_n$) such that, for $h_1 \geq  1$ and $h_j\geq  \varphi_{j-1}(h_{j-1})$, $j=2, \dots, n$ we have
	\[
	P\big(B_{n+1}=h_{n+1} \mid B_{n}=h_{n}, \dots, B_{1}=h_{1}\big)= F_n(\beta_n)-F_n(\alpha_n),
	\]
	where 
	\[
	\alpha_n=\delta_n(h_n, h_{n+1}+1, q_n)  ,\quad \beta_n=\delta_n(h_n, h_{n+1}, q_n)\quad\mbox{with}\quad\delta_j(h,k, q) = \frac{ \varphi_j (h )(1+q )}{k+\varphi_j (h ) q  }.
	\]
	Let $Q_n= q_n(B_1, \dots, B_n)$ and define
	\begin{equation}
		\label{Rdef}R_{n}= \frac{ B_{n+1}+\varphi_n(B_n) Q_n}{\varphi_n(B_n)(1+Q_n) }= \frac{1}{\delta_n(B_n, B_{n+1}, Q_n)}.
	\end{equation}

	\medskip 
	\noindent
	The origin of these random variables is related to digital expansions of real numbers: the sequence $(R_n)_{n \geq 1}$ describes the ratio of the digits and its definition unifies several particular cases that have been studied by different authors: the L\"uroth series was studied in \cite{G1976} and \cite{L1883}, the Engel series was studied in \cite{E1913} and \cite{S1974} while the cases of Engel continued fraction expansions and the Sylvester series were explored in \cite{HKS2002} and \cite{P1960} respectively. Several results for the general sequence $(R_n)_{n \geq 1}$ have been proven during the last decade:   in \cite{F2007} the Oppenheim continued fractions are introduced and their arithmetic and metric properties are studied, in \cite{G} a weak law has been established for the sequence $\displaystyle\frac{1}{n\log n}\sum_{k=1}^{n}R_k,$ while in \cite{GH} exact weak and strong laws (i.e. convergence either in probability or almost sure to a positive constant) have been proven for the sequence $\displaystyle \frac{1}{b_n}\sum_{k=1}^{n}a_kR_k$ for suitably chosen sequences of positive numbers $(a_n)_{n\geq 1}$ and $(b_n)_{n\geq 1}$. Moreover, in \cite{FG2021}, conditions are identified  under which $R_n^{-1}$ converges in distribution and the sequence $\displaystyle\sum_{k=1}^{n}\log R_k$ satisfies a central limit theorem.
	
	\medskip 
	\noindent	
	It is worth noting that for every integer $n$  the random variable $(R_n)_{n \geq 1}$ has a \textit{long-tailed} distribution. Long-tailed distributions are used in many different research areas such as in finance and actuarial, machine learning and artificial intelligence (see for example \cite{FKZ2011} for more on long-tailed distributions).
	
	\bigskip
	\noindent
	In this work we focus on studying the asymptotic behaviour of the extremes of generalized Oppenheim random variables, which is a research topic that has never being studied before in the literature.

	\medskip 
	\noindent
	In Sections 2 and 3 we provide results which will be used later for the proofs of the main theorems; these preliminary results include the asymptotic behaviour of some sequences arising from independent and identically distributed random variables (Lemma \ref{7}) and probability inequalities for the random variables $(R_n)_{n \geq 1}$ (Theorem \ref{6}). In Section 4 we study the asymptotic behaviour of extremes of a sequence of generalized Oppenheim random variables and in Section 5 we provide an Extremal Types Theorem for generalized Oppenheim expansions. 
	
	\medskip 
	\noindent
	In the whole paper by $(R_n)_{n \geq 1}$ we shall mean the sequence defined   in \eqref{Rdef} and we shall assume that the sequence $(F_n)_{n \geq 1}$ of the related distribution functions is constant for $n$, i.e., for every integer $n$, $F_n=F$ for some (not necessarily continuous) distribution function $F$ such that $$\lim_{t \to 0^+}F(t)=0 \quad {\rm and}  \quad \lim_{t \to 1^-}F(t)\leq 1.$$  
	Obviously this implies that $F(x)=0 $ for $x \leq 0$ and $F(x) =1$ for $x \geq 1$.  We shall also assume that the sequence of functions $(\varphi_n)_{n \geq 1}$ introduced above verifies the inequality $\varphi_n\geq 1$ for every integer $n$.  This assumption plays a crucial role in Lemma \ref{1} in the following section and also underpins several outcomes discussed in Section 5; see \cite{GH} for details.

	\medskip 
	\noindent
	For $x,y\in\mathbb{R}$ the notation $x\vee y$ stands for the $\max\{x,y\}$ while $I_{\{A\}}(x)$ will be used to denote the indicator function of the set $A$. 
	
	\section{Preliminaries}  
	We have the following multidimensional generalizations of Corollary 2.8 and Proposition 2.9 of   \cite{GH}:
	
	\begin{lemma}\label{1}\sl The following inequalities hold true: 
		\begin{enumerate}
			\item [(i)] for every $k$, every finite sequence of integers $i_1 , \dots,  i_k$  and  every finite sequence of numbers $x_{ 1}, \dots, x_{ k} \geq  1$ we have
			$$ F\left (\frac{1}{ x_{ 1}+1}\right ) \cdots F\left (\frac{1}{ x_{ k}+1}\right ) \leq P(R_{i_1}> x_{1}, \dots, R_{i_k}> x_{k})\leq F\left (\frac{1}{ x_{1}}\right ) \cdots F\left (\frac{1}{ x_{k}}\right). $$
			
			\item[(ii)] for every $k$, every finite sequence $I_{ 1} , \dots,  I_{ k}$  of intervals in $\mathbb{N}$   and  every finite sequence of numbers $x_{ 1}, \dots, x_{ k} \geq  1$ we have
			\begin{align*}&
				F^{q_1}\left (\frac{1}{ x_{1}+1}\right ) \cdots F^{q_k}\left (\frac{1}{ x_{k}+1}\right ) \leq P\left (\bigcap_{i \in I_{1}}\left \{R_{i }> x_{1}\right \}, \dots,  \bigcap_{i \in I_{k}}\left \{R_{i }> x_{k}\right \}\right )\\ &\leq F^{q_1}\left (\frac{1}{ x_{1}}\right ) \cdots F^{q_k}\left (\frac{1}{ x_{k}}\right ),
			\end{align*}			where, for every $j =1, \dots, k$, $q_j = \# I_j.$
		\end{enumerate}
	\end{lemma}
	\begin{proof}
		Point (ii) is a consequence of (i). The proof of (i) is by induction on $k$. The statement holds for $k=1$ and $k=2$ by Lemmas 2.5 and 2.6 in \cite{GH} respectively. Now assume that the statement holds for $k$; 
		by an argument similar to the one used in the proof of Lemma 2.5 of \cite{GH} and with the same notation  we have
		\begin{align}&\label{eq1}
			P(R_{i_1}> x_{i_1}, \dots, R_{i_k}> x_{i_{k+1}})\\& = \nonumber E\left[1_{\{R_{i_1}> x_{i_1}, \dots, R_{i_k}> x_{i_k}\}}F\left (\frac{\varphi_{i_{k+1}}(B_{i_{k+1}})(1+ Y_{i_{k+1}})}{S_{i_{k+1}}(x_{i_{k+1}};B_1, \dots, B_{i_{k+1}})+ \varphi_{i_{k+1}}(B_{i_{k+1}})Y_{i_{k+1}}}\right )\right],
		\end{align}where
		$S_n(x;B_1, \dots, B_n) $ is the random variable obtained by composition of the random vector $(B_1, \dots, B_{n})$  with the function $(h_1, \dots, h_{n})\mapsto s_n(x ;h_{1}, \dots, h_{n})$
		and
		$$s_n(x ;h_{1}, \dots, h_{n})=\lceil x \varphi_n (h_n)+ (x-1)y_n \varphi_n(h_n)\rceil.$$
		Under the given assumptions it is easy to see that
		$$\frac{1}{x+1}   \leq \frac{\varphi_{n}(B_{n})(1+ Y_{n})}{S_{n}(x;B_1, \dots, B_{n})+ \varphi_{n}(B_{n})Y_{n}}  \leq \frac{1}{x} $$
		so that, by \eqref{eq1},
		\begin{align*}&P(R_{i_1}> x_{i_1}, \dots, R_{i_k}> x_{i_{k+1}})\\&
			\leq 
			E\left[1_{\{R_{i_1}> x_{i_1}, \dots, R_{i_k}> x_{i_k}\}}\right]F\left (\frac{1}{x_{i_{k+1}}}\right )= P(R_{i_1}> x_{i_1}, \dots, R_{i_k}> x_{i_k})F\left (\frac{1}{x_{i_{k+1}}}\right )\\&\leq F\left (\frac{1}{ x_{i_1}}\right ) \cdots F\left (\frac{1}{ x_{i_k}}\right ) F\left (\frac{1}{x_{i_{k+1}}}\right ),
		\end{align*}
		by the inductive assumption. The left inequality runs similarly.
	\end{proof}

	\medskip
	\noindent 
	Let $F$ be the distribution function related to the sequence $(R_n)_{n\geq 1}$ (as explained in the Introduction) and let $(U_n)_{n\geq 1}$ be a sequence  of independent   and identically distributed random variables with common distribution $F$. In the result that follows we provide conditions under which some extremes arising from such sequence $(U_n)_{n\geq 1}$ belong to the MDA of the Weibull and the Fr\'echet distribution. This result will be used later in Section 4 for investigating the behaviour of the analogous extremes arising from $(R_n)_{n\geq 1}$.
	
	\begin{lemma}\label{7} \sl Consider a sequence of independent and identically distributed random variables $(U_n)_{n\geq 1}$ with common distribution $F$.
		
		\begin{enumerate}	
			\item[(i)]Assume the existence of the limit
			$$\lim_{t \to 1^-} \frac{F(t)-1}{t-1} =: \ell_{1 } ^-\in [0, + \infty].$$
			Then
			\begin{enumerate}
				\item if  $\ell_{1 } ^-\not= 0, + \infty$, the sequence $(U_n)_{n\geq 1}$ belongs to the MDA of  the negative Weibull distribution (with $\alpha =1$). More precisely $$\lim_{n \to \infty }    P\left (\max_{1 \leq k \leq n} U_k \leq 1 + \frac{x}{n}\right )=
				\begin{cases} e ^ {  \ell_{1 } ^- x}  & x \leq  0\\ 1 & x >0;
				\end{cases}
				$$
				
				\item if $\ell_{1 } ^-= + \infty$ we have $$\lim_{n \to \infty }    P\left (\max_{1 \leq k \leq n} U_k \leq 1 + \frac{x}{n}\right )=
				\begin{cases} 0  & x < 0\\ 1 & x \geq 0;
				\end{cases}
				$$
				
				\item if    $\ell_{1 } ^-= 0$ the above limit equals 1 for every $x \in \mathbb{R}$. 
			\end{enumerate}
			
			\item[(ii)] Assume the existence of the limit $$\lim_{t \to 0^+} \frac{F(t)}{t} =:  \ell_{0 } ^+\in [0, + \infty]. $$
			Then
			\begin{enumerate}
				\item if  $\ell_{0 } ^+\not= 0, + \infty$, the sequence $({U_n}^{-1})_{n\geq 1}$ belongs to the  MDA of  the Fr\'echet distribution (with $\alpha =1$). More precisely 
				$$\lim_{n \to \infty }P\left (\max_{1 \leq k \leq n} \frac{1}{U_k} \leq   xn\right )   = \begin{cases} 0  & x\leq 0\\ e ^ { -\frac{\ell_{0 }^+}{x} }  & x > 0 ;
				\end{cases} 
				$$
				
				\item if $ \ell_{0 } ^+= 0$ we have $$\lim_{n \to \infty }    P\left (\max_{1 \leq k \leq n}  \frac{1}{U_k} \leq xn\right )=
				\begin{cases} 0  & x \leq 0\\ 1 & x > 0;
				\end{cases}
				$$
				
				\item if  $ \ell_{0 } ^+= +\infty$ the above limit equals 0 for every $x \in \mathbb{R}$. 
			\end{enumerate}
		\end{enumerate}
	\end{lemma}
	\begin{proof}
		In (i) (a) and (ii) (a) the fact that the two sequences belong to the MDA of the Weibull and the Fr\'echet distribution respectively is a consequence of the Fisher-Tippett-Gnedenko Theorem. Anyway we are interested in the explicit calculations, which work also in the other cases (i.e. (b) and (c) of both (i) and (ii),  the proofs of  which are omitted).
		
		\bigskip
		\noindent
		For (i) (a) the  case $ x \geq 0$ is obvious. For $x < 0$ we have
		\begin{align*}&
			P\left (\max_{1 \leq k \leq n} U_k \leq 1 + \frac{x}{\ell_{1 } ^-n}\right )= P\left (\bigcap_{1 \leq k \leq n}\left \{U_k \leq 1 + \frac{x}{\ell_{1 } ^-n}\right \}\right )
			= \prod_{1 \leq k \leq n} P\left (U_k \leq 1 + \frac{x}{\ell_{1 } ^-n}\right ) \\&  =F^n\left (1 + \frac{x}{\ell_{1 } ^-n}\right ).\end{align*}  
		Now
		\begin{align*}& n\log F \left (1 + \frac{x}{\ell_{1 } ^-n}\right ) = n\log \left (1+\left \{ F \left (1 + \frac{x}{\ell_{1 } ^-n}\right )-1\right \}\right )\sim n\left( F \left (1 + \frac{x}{\ell_{1 } ^-n}\right )-1\right)\\ & = \frac{F \left (1 + \frac{x}{\ell_{1 } ^-n}\right )-1}{\frac{x}{\ell_{1 } ^-n}}\cdot \frac{x}{ \ell_{1 } ^-}\to   x, \qquad n \to \infty.
		\end{align*} 
		
		\noindent
		For (ii) (a) the  case $ x \leq 0$ is obvious. For $x > 0$,  $xn\ell_{0 }^+$ is ultimately larger than 1;  thus, ultimately, 
		\begin{align*}&
			P\left (\max_{1 \leq k \leq n} \frac{1}{U_k} \leq   xn\ell_{0 }^+\right )   =P\left (\min_{1 \leq k \leq n} U_k  \geq \frac{1}{xn\ell_{0 }^+}\right ) =P\left (\bigcap_{1 \leq k \leq n}\left \{U_k \geq  \frac{1}{xn\ell_{0 }^+} \right \}\right )= \left (1 - F \left (\frac{1}{xn\ell_{0 }^+} \right ) \right )^n.\end{align*} 
		Now
		\begin{align*}& n \log \left \{1 - F \left (\frac{1}{xn\ell_{0 }^+} \right ) \right \}\sim - n F \left (\frac{1}{xn\ell_{0 }^+} \right )  = -  \frac{F (\frac{1}{xn}\ell_{0 }^+ ) }{\frac{1}{xn\ell_{0 }^+}} \cdot\frac{1}{x\ell_{0 }^+} \to -\frac{1}{x}, \qquad n \to \infty.
		\end{align*}  
	\end{proof} 
	
	\begin{remark}
		It is worth noting that the result in (i) (a) can be rephrased to
		$$\lim_{n \to \infty }    P\left (\max_{1 \leq k \leq n} U_k \leq 1 + \frac{x}{ \ell_{1 } ^-n}\right )=
		\begin{cases} e ^ {   x}  & x \leq 0\\ 1 & x > 0;
		\end{cases}$$
		while the expression in (ii) (a) can be equivalently written as
		$$\lim_{n \to \infty }P\left (\max_{1 \leq k \leq n} \frac{1}{U_k} \leq   x n\ell_{0 }^+ \right )   = \begin{cases} 0  & x\leq 0\\ e ^ {   -  \frac{1}{x}}  & x > 0.
		\end{cases} $$
		
	\end{remark}

	\section{Some useful inequalities}
	Let $a\geq 1$. In the following  result we use the notation
	$$S_a = \displaystyle\frac{ F (\frac{1}{a})- F (\frac{1}{a+1}) }{ \frac{1}{a}-\frac{1}{a+1}}.$$
	
	\begin{theorem}\label{6}\sl Let  $a, b$  be  constants.
		\begin{enumerate}
			\item[(i)] For $1 \leq a < b$  we have,
			\begin{align*}&
				\left|P\left(\bigcap_{1 \leq k \leq n}\left\{a < R_k \leq b\right\}\right)- F^n\left (\frac{1}{a}\right ) \left( 1-F\left (\frac{1}{b}\right )\right)^n\right|  \leq \alpha^{(a,b)}_n,
			\end{align*}
			where
			\begin{align*}
				\alpha^{(a,b)}_n&=  \frac{(S_a \vee S_b) n}{b^2}\left(1 + F\left (\frac{1}{b}\right )\right)^{n-1} F^{n}\left (\frac{1}{a}\right ) +\frac{ (S_a \vee S_b) n}{a^2}  F^{n-1}\left (\frac{1}{a}\right )\left(1 + F\left (\frac{1}{b}\right )\right)^{n}\\
				& +\frac{(S_a \vee S_b) n}{a^2}   F^{n-1}\left (\frac{1}{a}\right ) ,
			\end{align*}  
			
			\item[(ii)] Let $b\geq 1$. Then
			$$\left|P\left (\bigcap_{1 \leq k \leq n}\left \{R_k \leq b\right \}\right )-\left(1-F\left (\frac{1}{b}\right )\right)^n \right|\leq  \frac{ S_b   n}{b^2}\left(1 + F\left (\frac{1}{b}\right )\right)^{n-1}.$$

			\item[(iii)]  Let $a\geq 1$. Then 
			$$\left|P\left (\bigcap_{1 \leq k \leq n}\left \{R_k > a\right \}\right )- F^n\left (\frac{1}{a}\right )  \right|\leq  \frac{ S_a  n}{a^2}  F^{n-1}\left (\frac{1}{a}\right ).$$
		\end{enumerate}
	\end{theorem}
	\begin{proof}
		We provide the proof only for part (i) as (ii) and (iii) can be obtained in a similar manner.
		First,  denote $$A_k = \{R_k > a\}, \quad  B_k = \{R_k > b\}, \quad k = 1, \dots, n, $$ and  $$A=\bigcap_{1 \leq k \leq n}A_k, \quad B=  \bigcup_{1 \leq k \leq n} B  _k.$$
		Then the event $\displaystyle\bigcap_{1 \leq k \leq n}\{a < R_k \leq b\}$
		can be written in the form 
		$$\bigcap_{1 \leq k \leq n}\left (A_k\cap B^c _k \right )= A \cap \left (\bigcap_{1 \leq k \leq n} B^c _k\right )= A \cap B^c.$$
		Thus, by the inclusion-exclusion formula,
		\begin{align*}&
			P\left (\bigcap_{1 \leq k \leq n}\left \{a < R_k \leq b\right \}\right)= P(A) - P(A \cap B)= P(A)- P\left (\bigcup_{k=1}^n(A \cap B_k)   \right )\\& = P(A)- \left \{\sum_{i=1}^n (-1)^{i+1} \sum_{1 \leq j_1 < \cdots < j_i\leq n} P\left (A \cap \bigcap _{r=1}^i B_{j_r} \right )\right \}\\& =P(A)+ \sum_{1 \leq i \leq n\atop even \, i}\sum_{1 \leq j_1 < \cdots < j_i\leq n} P\left (A \cap\bigcap _{r=1}^i B_{j_r}\right )-  \sum_{1 \leq i \leq n\atop odd \, i}\sum_{1 \leq j_1 < \cdots < j_i\leq n} P\left (A \cap\bigcap _{r=1}^i B_{j_r}\right )  \end{align*}
		By Lemma \ref{1} we obtain that
		\begin{align*}&
			P\left (\bigcap_{1 \leq k \leq n}\left \{a < R_k \leq b\right \}\right )\\&\leq F^n\left (\frac{1}{a}\right )+ \sum_{1 \leq i \leq n\atop even \, i}\sum_{1 \leq j_1 < \cdots < j_i\leq n}F^n\left (\frac{1}{a}\right )F^i\left (\frac{1}{b}\right ) -  \sum_{1 \leq i \leq n\atop odd \, i}\sum_{1 \leq j_1 < \cdots < j_i\leq n} F^n\left (\frac{1}{a+1}\right )F^i\left (\frac{1}{b+1}\right ) \\& = \underbrace{\left \{F^n\left (\frac{1}{a}\right )+ \sum_{1 \leq i \leq n\atop even \, i}\sum_{1 \leq j_1 < \cdots < j_i\leq n}F^n\left (\frac{1}{a}\right )F^i\left (\frac{1}{b}\right ) -  \sum_{1 \leq i \leq n\atop odd \, i}\sum_{1 \leq j_1 < \cdots < j_i\leq n} F^n\left (\frac{1}{a}\right )F^i\left (\frac{1}{b}\right )\right \}}_{= C_n }\\& + \underbrace{ \left \{\sum_{1 \leq i \leq n\atop odd \, i}\sum_{1 \leq j_1 < \cdots < j_i\leq n} \left (F^n\left (\frac{1}{a}\right )F^i\left (\frac{1}{b}\right )-F^n\left (\frac{1}{a+1}\right )F^i\left (\frac{1}{b+1}\right )\right )\right \}}_{= D_n }.
		\end{align*}
		Now
		\begin{align*}& 
			C_n= F^n\left (\frac{1}{a}\right ) + \sum_{i=1}^n (-1)^i\sum_{1 \leq j_1 < \cdots < j_i\leq n}F^n\left (\frac{1}{a}\right )F^i\left (\frac{1}{b}\right )\\&= F^n\left (\frac{1}{a}\right ) +\sum_{i=1}^n (-1)^iF^n\left (\frac{1}{a}\right )F^i\left (\frac{1}{b}\right )\sum_{1 \leq j_1 < \cdots < j_i\leq n}1\\&= F^n\left (\frac{1}{a}\right ) +\sum_{i=1}^n (-1)^i{n \choose i}F^n\left (\frac{1}{a}\right )F^i\left (\frac{1}{b}\right )= \sum_{i=0}^n (-1)^i{n \choose i}F^n\left (\frac{1}{a}\right )F^i\left (\frac{1}{b}\right )\\&= F^n\left (\frac{1}{a}\right )\left(1-F\left (\frac{1}{b}\right )\right)^n.
		\end{align*}
		For bounding the term $D_n$, we first study the term 
		\[
		F^n\left (\frac{1}{a}\right )F^i\left (\frac{1}{b}\right )-F^n\left (\frac{1}{a+1}\right )F^i\left (\frac{1}{b+1}\right ).
		\]
		Observe that
		\begin{align}
			\nonumber &F^n\left (\frac{1}{a}\right )F^i\left (\frac{1}{b}\right )-F^n\left (\frac{1}{a+1}\right )F^i\left (\frac{1}{b+1}\right )\\
			\nonumber 
			& = F^n\left (\frac{1}{a}\right )\left (F^i\left (\frac{1}{b}\right ) - F^i\left (\frac{1}{b+1}\right )\right )+F^i\left (\frac{1}{b+1}\right )\left (F^n\left (\frac{1}{a}\right ) - F^n\left (\frac{1}{a+1}\right )\right )\\
			\nonumber 
			&\leq F^n\left (\frac{1}{a}\right )\left (F^i\left (\frac{1}{b}\right ) - F^i\left (\frac{1}{b+1}\right )\right )+F^i\left (\frac{1}{b}\right )\left (F^n\left (\frac{1}{a}\right ) - F^n\left (\frac{1}{a+1}\right )\right )\\
			\label{3}
			&\leq nF^{n-1}\left (\frac{1}{a}\right )F^i\left (\frac{1}{b}\right )\left(F\left (\frac{1}{a}\right ) - F\left (\frac{1}{a+1}\right )\right)+iF^n\left (\frac{1}{a}\right )F^{i-1}\left (\frac{1}{b}\right )\left(F\left (\frac{1}{b}\right ) - F\left (\frac{1}{b+1}\right )\right)
		\end{align}
		where \eqref{3} follows from the inequality
		\[
		0<x^n - y^n= (x-y)(x^{n-1}+x^{n-2}y+\cdots+y^{n-1}) \leq nx^{n-1}(x-y),\quad\mbox{where}\quad x>y
		\]
		for integer values of $n\geq 1$.
		The above is  less or equal to
		\begin{align} \label{5}&  \frac{ (S_a \vee S_b) n}{a^2}   F^{n-1}\left (\frac{1}{a}\right )  F^i\left (\frac{1}{b}\right )+ \frac{(S_a \vee S_b) i}{b^2} F^{n}\left (\frac{1}{a}\right )  F^{i-1}\left (\frac{1}{b}\right );
		\end{align}
		hence
		\begin{align}& \label{4}
			D_n\leq  (S_a \vee S_b) \sum_{1 \leq i \leq n\atop odd \, i}\sum_{1 \leq j_1 < \cdots < j_i\leq n}  \left \{  \frac{ n}{a^2}   F^{n-1}\left (\frac{1}{a}\right )  F^i\left (\frac{1}{b}\right )+ \frac{ i}{b^2} F^{n}\left (\frac{1}{a}\right )  F^{i-1}\left (\frac{1}{b}\right )\right \}\\
			& = \nonumber  \frac{ S_a \vee S_b }{a^2}  F^{n-1}\left (\frac{1}{a}\right ) \sum_{1 \leq i \leq n\atop odd \, i}{n \choose i} F^i\left (\frac{1}{b}\right )+ \frac{S_a \vee S_b}{b^2} F^{n}\left (\frac{1}{a}\right ) \sum_{1 \leq i \leq n\atop odd \, i}{n \choose i} i  F^{i-1}\left (\frac{1}{b}\right )\\
			&\leq   \nonumber
			\frac{S_a \vee S_b}{a^2}  F^{n-1}\left (\frac{1}{a}\right ) \sum_{i=0}^n{n \choose i} F^i\left (\frac{1}{b}\right )+ \frac{S_a \vee S_b}{b^2} F^{n}\left (\frac{1}{a}\right )  \sum_{i=1}^n {n \choose i} i  F^{i-1}\left (\frac{1}{b}\right )
			\\& =   \nonumber \frac{(S_a \vee S_b )n}{a^2}  F^{n-1}\left (\frac{1}{a}\right )\left \{1+F\left (\frac{1}{b}\right ) \right \}^n+ \frac{(S_a \vee S_b )n}{b^2}  F^{n}\left (\frac{1}{a}\right )\left \{1+F\left (\frac{1}{b}\right ) \right \}^{n-1}.
		\end{align}
		since, for every $z \in \mathbb{R}$, \begin{align*}&
			\sum_{i=1}^n{n \choose i} iz^{i-1}= n\sum_{i=1}^n \frac{(n-1)!}{(i-1)!(n-i)!}z^{i-1}= n\sum_{i=0}^{n-1}{n-1 \choose i}z^i = n(1+z)^{n-1}.
		\end{align*}
		The LHS inequality is obtained by similar calculations. Precisely
		\begin{align*}&
			P\left(\bigcap_{1 \leq k \leq n}\left\{a < R_k \leq b\right\}\right)\\
			& \geq F^n\left (\frac{1}{a+1}\right )+ \sum_{1 \leq i \leq n\atop even \, i}\sum_{1 \leq j_1 < \cdots < j_i\leq n}F^n\left (\frac{1}{a+1}\right )F^i\left (\frac{1}{b+1}\right ) -  \sum_{1 \leq i \leq n\atop odd \, i}\sum_{1 \leq j_1 < \cdots < j_i\leq n} F^n\left (\frac{1}{a}\right )F^i\left (\frac{1}{b}\right ) \\
			& =\underbrace{ F^n\left (\frac{1}{a+1}\right )- F^n\left (\frac{1}{a}\right )}_{=E_n} \\&+\underbrace{\left \{F^n\left (\frac{1}{a}\right )+ \sum_{1 \leq i \leq n\atop even \, i}\sum_{1 \leq j_1 < \cdots < j_i\leq n}F^n\left (\frac{1}{a}\right )F^i\left (\frac{1}{b}\right ) -  \sum_{1 \leq i \leq n\atop odd \, i}\sum_{1 \leq j_1 < \cdots < j_i\leq n} F^n\left (\frac{1}{a}\right )F^i\left (\frac{1}{b}\right )\right \}}_{= C_n }\\
			& + \underbrace{ \left \{\sum_{1 \leq i \leq n\atop even \, i}\sum_{1 \leq j_1 < \cdots < j_i\leq n} \left (F^n\left (\frac{1}{a+1}\right )F^i\left (\frac{1}{b+1}\right )-F^n\left (\frac{1}{a}\right )F^i\left (\frac{1}{b}\right )\right )\right \}}_{= H_n }.
		\end{align*}
		Now (as before)
		\begin{align*}& 
			C_n=   F^n\left(\frac{1}{a}\right)\left(1-F\left(\frac{1}{b}\right)\right)^n;
		\end{align*}
		and (see \eqref{5})
		\begin{align*}&F^n\left (\frac{1}{a}\right )F^i\left (\frac{1}{b}\right )-F^n\left (\frac{1}{a+1}\right )F^i\left (\frac{1}{b+1}\right ) \\& \leq \frac{(S_a \vee S_b) n}{a^2}   F^{n-1}\left (\frac{1}{a}\right )  F^i\left (\frac{1}{b}\right )+ \frac{(S_a \vee S_b)i}{b^2} F^{n}\left (\frac{1}{a}\right )  F^{i-1}\left (\frac{1}{b}\right ).
		\end{align*}
		hence (similarly as above, see \eqref{4})
		\begin{align*}& 
			- H_n\leq   (S_a \vee S_b)\sum_{1 \leq i \leq n\atop even \, i}\sum_{1 \leq j_1 < \cdots < j_i\leq n} \left \{  \frac{n}{a^2}   F^{n-1}\left (\frac{1}{a}\right )  F^i\left (\frac{1}{b}\right )+ \frac{i}{b^2} F^{n}\left (\frac{1}{a}\right )  F^{i-1}\left (\frac{1}{b}\right )\right \}\\
			& =   \frac{(S_a \vee S_b) n}{a^2}  F^{n-1}\left (\frac{1}{a}\right ) \sum_{1 \leq i \leq n\atop even \, i}{n \choose i} F^i\left (\frac{1}{b}\right )+ \frac{S_a \vee S_b}{b^2} F^{n}\left (\frac{1}{a}\right ) \sum_{1 \leq i \leq n\atop even \, i}{n \choose i} i  F^{i-1}\left (\frac{1}{b}\right )\\&\leq  
			\frac{(S_a \vee S_b) n}{a^2}  F^{n-1}\left (\frac{1}{a}\right ) \sum_{i=0}^n{n \choose i} F^i\left (\frac{1}{b}\right )+ \frac{S_a \vee S_b}{b^2} F^{n}\left (\frac{1}{a}\right )  \sum_{i=1}^n {n \choose i} i  F^{i-1}\left (\frac{1}{b}\right )
			\\& =   \frac{(S_a \vee S_b) n}{a^2}  F^{n-1}\left (\frac{1}{a}\right )\left \{1+F\left (\frac{1}{b}\right ) \right \}^n+ \frac{(S_a \vee S_b) n}{b^2}  F^{n}\left (\frac{1}{a}\right )\left \{1+F\left (\frac{1}{b}\right ) \right \}^{n-1}.\end{align*}
		The case $i=0$ of \eqref{3} reads as
		$$- E_n =F^n\left (\frac{1}{a}\right ) -F^n\left (\frac{1}{a+1}\right )  \leq \frac{(S_a \vee S_b) n}{a^2}   F^{n-1}\left (\frac{1}{a}\right )   $$
		and the LHS of the inequality is proved.
		
	\end{proof}
	
	\section{The asymptotic behaviour of some extremes} 
	Assume that there exists $\beta > 0$  such that 
	\begin{equation}
		\label{condF}F(x) - F(y) \leq \beta (x-y)\quad \mbox{for}\quad x\geq y.
	\end{equation}
	It is worth mentioning that this assumption enables us to circumvent the need for the highly stringent assumption of the differentiability of the involved distribution function.
	
	\medskip
	\noindent Theorem \ref{6} has some interesting consequences concerning the sequences $$M_n :=\max_{1 \leq k \leq n}  R_k, \qquad Z_n  :=\min_{1 \leq k \leq n}  R_k$$
	which are discussed in the results of this section.
	
	\begin{corollary}\label{14}\sl Assume that $\ell_0^+$ exists as a finite number   and condition \eqref{condF} holds true. Then 
		\begin{enumerate}
			\item [(i)] if $\ell_0^+>0$, the sequence $(M_n )_n$  belongs to the   MDA of the   Fr\'echet distribution (with $\alpha =1$). Precisely 
			$$\lim_{n \to \infty }P\left ( M_n  \leq   xn\ell_0^+\right )   = \begin{cases} 0  & x\leq0\\ e ^ {  -\frac{1}{x}}  & x > 0 ;
			\end{cases} 
			$$
			\item[(ii)] if $\ell_0^+=0$, then
			$$\lim_{n \to \infty }P\left ( M_n  \leq   xn\right )   = \begin{cases} 0  & x\leq0\\1  & x > 0 .
			\end{cases} 
			$$
		\end{enumerate}
	\end{corollary} 
	\begin{proof}The case $x \leq 0$ is obvious due to the fact that $R_n \geq 1$ for every $n$. For every $x > 0$, $xn \ell_0^+ $ is ultimately  larger than 1; thus, by (ii) of the Theorem \ref{6},
		\begin{align*}&
			\left|P\left (\bigcap_{1 \leq k \leq n}\left \{R_k \leq xn\ell_0^+\right \}\right )-\left(1-F\left (\frac{1}{xn\ell_0^+}\right )\right)^n \right| \leq \frac{\beta }{nx^2(\ell_0^+)^2}\left(1 + F\left (\frac{1}{nx\ell_0^+}\right )\right)^{n-1} \to 0
		\end{align*}
		by an argument similar to the one used in the proof of  Lemma \ref{7}, point (ii) (see last line  there).
		Since
		$$\left( 1-F\left (\frac{1}{xn\ell_0^+}\right )\right)^n= P\left (\min_ {1 \leq k \leq n}U_k > \frac{1}{xn\ell_0^+}\right )= P\left (\max _{1 \leq k \leq n}\frac{1}{U_k}< n \ell_0^+x\right ),$$
		we conclude by Lemma \ref{7}, point (ii).
	\end{proof}
	
	\begin{remark} This corollary implies that if the sequence $(R_n)_n$ is uniformly bounded by a constant, then $\ell_0^+=0$.
	\end{remark}
	
	\begin{corollary}\label{11}\sl  Assume  that the distribution function $F$ satisfies condition \eqref{condF} and  
		
		\begin{equation}\label{eq8}
			\liminf_{y \downarrow 0} \frac{\log F(1- h y)}{y \log y}> 1, \qquad \, \forall h> 0.
		\end{equation}
		Then  
		\begin{enumerate}
			\item[(i)] $\ell_1^- =+\infty;$
			
			\item[(ii)]we have
			$$\lim_{n \to \infty }P\left (  \frac{1}{Z_n}  \leq  1+\frac{x}{n}\right )   = \begin{cases} 0    & x< 0\\  1 & x \geq 0. 
			\end{cases} 
			$$
		\end{enumerate}
	\end{corollary} 
	\begin{proof}
		For  part (i) we start by assuming the contrary i.e. that $\ell_1^-\in \mathbb{R}^+$. Fix any $h >0$; then, as  $y \downarrow 0$,
		\begin{align}& \frac{\log F(1- h y)}{y \log y} \sim \frac{  F(1- h y)- 1}{y \log y}= \frac{F(1- h y)- 1}{-hy} \cdot\frac{-h}{ \log y}\to \ell_1^-\cdot 0 =0,
		\end{align}
		which contradicts assumption   \eqref{eq8}.  
		
		\noindent Next, for part (ii), the case $x  \geq 0$ is obvious.  Let $x < 0$. Notice that $$\frac{1}{Z_n }=\max_{1 \leq k \leq n} \frac{1}{R_k}.$$ Hence, by (iii) of  Theorem \ref{6},
		\begin{align*}& 
			\left|P\left (\bigcap_{1 \leq k \leq n}\left \{\frac{1}{R_k}  \leq 1+\frac{x}{n} \right \}\right )- F^n\left (1+\frac{x}{n}\right )  \right|\leq  \beta n    F^{n-1}\left (1+\frac{x}{n}\right ) 
		\end{align*}
		and we  prove that $$n   F^{n-1}\left (1+\frac{x}{n}\right )  \to 0,\qquad n \to \infty,$$  
		or equivalently that
		$$\log n + (n-1) \log F\left (1 + \frac{x}{n}\right ) \to - \infty.$$
		This is true since the above can be put in the form
		$$\log n\left \{1 + (n-1)\frac{\log F(1 + \frac{x}{n})}{\log n}  \right \},$$
		which is ultimately less or equal to  $- \epsilon \log n$ for some $\epsilon > 0, $ by assumption   \eqref{eq8}. Since 
		$$F^n\left (1+\frac{x}{n}\right )= P\left (\max_{1 \leq k \leq n}  U_k \leq 1+\frac{x}{n}\right ),$$we conclude by Lemma \ref{5}, point (i)(b).
		
	\end{proof}

	\begin{remark}  Obviously the assumption \eqref{eq8} can be equivalently rephrased by requiring  that there exists $ \tilde h >0$ such that the inequality holds for every $0<h< \tilde h$.
	\end{remark}
	
	\medskip
	
	\begin{remark} Assumption \eqref{eq8} is somewhat unpleasant since it depends on the parameter $h$; thus we present a condition assuring \eqref{eq8} and depending only on $F$.  Particularly, assume that $F$ is  differentiable with derivative $f$ on $(0,1)$; assume moreover that $f$ is increasing and that
		\begin{equation}\label{eq10}
			\lim_{y \downarrow 0} \frac{f(1-y)}{-\log y} = + \infty.
		\end{equation}
		Then, assumption  \eqref{eq8} holds true. To prove this, fix any $h >0$; then, by Lagrange Theorem, there exists $\xi(y)$ with $1- hy < \xi(y) < 1$ such that
		\begin{align}&  \label{eq9} 
			\frac{\log F(1- h y)}{y \log y} = \frac{\log F(1- h y)- \log F(1)}{y \log y}= \frac{-hy f(\xi(y))}{F(\xi(y))y \log y} =  \frac{-h  f(\xi(y))}{F(\xi (y))  \log y}\nonumber\\&\sim \frac{-h f(\xi(y))}{ \log y} , \qquad y \downarrow 0
		\end{align}
		(the equivalence  $\sim $ holds since $F(1-hy )\leq  F(\xi(y))\leq  F(1)=1$). Since
		
		$$ \liminf_{y \downarrow 0}\frac{-h f(\xi(y))}{ \log y}\geq \liminf_{y \downarrow 0}\frac{-h f( 1-hy)}{ \log y},$$
		it suffices to ensure that
		$$\liminf_{y \downarrow 0}\frac{-h f( 1-hy)}{ \log y}> 1.$$
		
		\noindent Observe that, for sufficiently small $y$ we have
		$$\frac{-h f( 1-hy)}{ \log y} = \frac{-h f( 1-hy)}{ \log (hy)}\cdot \dfrac{\log(hy)}{\log y} =\frac{-h f( 1-hy)}{ \log (hy)} \left (1+\dfrac{\log h}{\log y}\right ) > \frac{2}{h}\cdot  h= 2,$$
		where the inequality is due to assumption \eqref{eq10} and the fact that $\frac{\log h}{\log y} \to 0$ as $y \downarrow 0$.
		
		\medskip
		
		\noindent Note that we can easily identify distributions that satisfy the proposed condition. For example, let $\alpha \in (0, 1)$ and take $$F(x) = 1 - (1-x)^ \alpha, \qquad 0 \leq x < 1.$$
		Then $f(x) = \frac{\alpha}{(1-x)^{1-\alpha}} $ is increasing and assumption  \eqref{eq10} is satisfied.
	\end{remark}
	
	\bigskip
	\bigskip

	\begin{remark}\label{15} The following example shows that   Corollary \ref{11} may be false without assuming \eqref{eq8}.  
		Let $(R_n)_{n \geq 1}$ be the L\"uroth sequence, i.e. $(R_n)_{n \geq 1}$ are independent with density
		$$P(R_n=h)= \frac{1}{h(h-1)}, \qquad h=2, 3, \dots.$$
		As it is well known (see \cite{G}), in this case the related distribution function is the uniform distribution on $[0,1]$, which doesn't satisfy \eqref{eq8}. For $t\geq 1 $ we have 
		$$  P(R_n \geq  t) = \sum_{h=\lceil t\rceil}^{\infty} \frac{1}{h(h-1)} = \sum_{h=\lceil t\rceil }^{\infty} \left (\frac{1}{h -1} - \frac{1}{h }\right )= \dfrac{1}{\lceil t\rceil-1}, $$ 
		where $\lceil t\rceil$ represents the smallest integer greater than or equal to $t$. 
		Let $p$ be an integer.	Then, for $x < 0$ and sufficiently large $n$ (depending on $x$ only)
		$$ P\left (\frac{1}{Z_n}  \leq  \frac{1}{p}+\frac{x}{n}\right) =P\left(\min_{1\leq k \leq n}R_k \geq \frac{np}{n+xp}\right)= \left(\frac{1}{\lceil \frac{np}{n+xp}\rceil-1 }\right)^n = \begin{cases}1 &   \hbox{ if} \,\,p=1\\ 0  &  \hbox{ if} \,\,p>1.
		\end{cases}$$
		This example shows that the behaviour of $(  Z_n^{-1})_{n \geq 1}$ may be heavily affected by the   normalization. Notice also that for no values of $p$ we have  found the behaviour described in Lemma \ref{7}, point (i) (a). 
	\end{remark}	
	
	\noindent 
	In what follows  we provide two results that establish some kind of \lq\lq asymptotic independence\rq\rq \ between the two sequences $(M_n)_{n\geq 1}$ and $(Z_n)_{n\geq 1}$.

	\begin{corollary}\label{12}\sl Assume that $1<x<y$ and  that the distribution function $F$ satisfies condition \eqref{condF} and   $F\left (\frac{1}{x}\right )< \frac{1}{2}$. Then 
		$$\lim_{n \to \infty}\left \{P(Z_n>x , M_n \leq y )-P(Z_n>x )P(M_n \leq y)\right \}=0.$$
	\end{corollary} 
	\begin{proof}
		Denote for simplicity
		\begin{align*}&
			a_n= P(Z_n>x),\quad  b_n = P(M_n \leq y), \quad c_n= P(Z_n>x , M_n \leq y ), \\&   \alpha_n =F^n\left (\frac{1}{x} \right ) ,\beta_n =\left ( 1-F\left (\frac{1}{y}\right )\right )^n  ,
		\end{align*}
		and write
		$$|c_n - a_n b_n |\leq  |c_n - \alpha_n \beta_n|+\alpha_n | b_n-\beta_n | +b_n|a_n - \alpha_n| $$
		and we bound the three summands by   means of Theorem \ref{6}. Particularly, for the first summand  we have that
		\begin{align*}&
			|c_n - \alpha_n \beta_n| \leq \frac{(S_x\vee S_y)n}{y^2}   F^n\left (\frac{1}{x} \right )\left (1+F\left (\frac{1}{y}\right )\right )^{n-1}\\&+  \frac{(S_x\vee S_y)n}{x^2}
			F^{n-1}\left (\frac{1}{x} \right )\left ( 1+F\left (\frac{1}{x}\right )\right )^{n }+ \frac{(S_x\vee S_y)n}{x^2} F^{n-1}\left (\frac{1}{x} \right )\\
			&\leq 
			\frac{(S_x\vee S_y) n}{ {x}^2}\left (1+F\left (\frac{1}{ {x}}\right )\right )^{n-1} F^{n-1}\left (\frac{1}{x} \right ) \left ( 2 F\left (\frac{1}{x} \right )+1\right )+ \frac{(S_x\vee S_y)n}{x^2} F^{n-1}\left (\frac{1}{x} \right )\\
			&\leq{ \frac{2(S_x\vee S_y) n}{x^2}\left (\dfrac{1}{2}\right )^{n-1}\left (\dfrac{3}{2}\right )^{n-1} +\frac{(S_x\vee S_y)n}{2^{n-1}x^2}\to 0,\quad n\to \infty,} 
		\end{align*}
		where the last two inequalities are due to the two assumptions on $x$.  For the second term and by employing again the assumptions on $x$  we have 
		\begin{align*}&
			\alpha_n | b_n-\beta_n |\leq \frac{S_y n}{y^2} F^{n}\left (\frac{1}{x} \right )\left ( 1+F\left (\frac{1}{y}\right )\right )^{n-1}\leq \frac{S_y n}{y^2}\left (\dfrac{1}{2}\right )^n \left (\dfrac{3}{2}\right )^{n-1} \to 0,
		\end{align*}  
		while similar arguments confirm the convergence to zero of the last summand 
		\begin{align*}& b_n|a_n - \alpha_n|\leq |a_n - \alpha_n|\leq  \frac{S_x n}{x^2}\left (F \left (\frac{1}{x}\right  )\right )^{n-1}\to 0.\end{align*} 
		
	\end{proof}
	
	\begin{remark} 
		If $\lim_{x \to 1^-}F(x) < \frac{1}{2}$,   the statement of   Corollary \ref{12} holds for every $1<x< y$, i.e. $(M_n)_n$ and $(Z_n)_n$ are asymptotically independent.
	\end{remark}
	
	\bigskip
	
	\noindent The following result can be proved similarly to Corollary \ref{12}.
	
	\begin{corollary}\sl  Assume that   condition \eqref{condF}  is satisfied and   there  exists $A>1  $ such that
		\begin{equation}\label{10}
			\sup_{a\geq A }S_a = C< \infty.
		\end{equation}
		Then $\dfrac{M_n}{\sigma_n}$ and $\dfrac{Z_n}{\rho_n}$ are asymptotically independent for every pair of sequences $(\rho_n)_n$ and $(\sigma_n)_n$ such that
		$$\lim_{n \to \infty }\rho_n =   \lim_{n \to \infty }\sigma_n = \infty$$ i.e.
		$$\lim_{n \to \infty}\left \{P(Z_n>x\rho_n, M_n \leq y\sigma_n)-P(Z_n>x \rho_n)P(M_n \leq y\sigma_n)\right \}=0 $$ 
	\end{corollary} 
	\begin{proof}
		Denote for simplicity
		\begin{align*}&
			a_n= P(Z_n>x\rho_n),\quad  b_n = P(M_n \leq y\sigma_n), \quad c_n= P(Z_n>x\rho_n , M_n \leq y\sigma_n ), \\&   \alpha_n =F^n\left (\frac{1}{x\rho_n} \right ) ,\beta_n =\left( 1-F\left (\frac{1}{y\sigma_n}\right )\right)^n  .
		\end{align*}
		Write
		$$|c_n - a_n b_n |\leq  |c_n - \alpha_n \beta_n|+\alpha_n |b_n-\beta_n|  + b_n|a_n -\alpha_n|.$$
		We bound the three summands by means of Theorem \ref{6}. For the first term, 
		\begin{align*}&
			|c_n - \alpha_n \beta_n| \leq \frac{(S_{x\rho_n}\vee S_{y\sigma_n})n}{(y \sigma_n)^2}   F^n\left (\frac{1}{x\rho_n} \right )\left ( 1+F\left (\frac{1}{{y\sigma_n}}\right )\right )^{n-1}\\&+  \frac{(S_{x\rho_n}\vee S_{y\sigma_n})n}{(x\rho_n)^2}
			F^{n-1}\left (\frac{1}{x\rho_n} \right )\left ( 1+F\left (\frac{1}{x\rho_n}\right )\right )^{n }+ \frac{(S_{x\rho_n}\vee S_{y\sigma_n})n}{(x\rho_n)^2} F^{n-1}\left (\frac{1}{x\rho_n} \right )\\
			& \leq \frac{ Cn}{(y \sigma_n)^2} \left (2F\left (\frac{1}{x\rho_n} \right ) \right )^{n-1}+  \frac{ Cn}{(x \rho_n)^2} \left (2F\left (\frac{1}{x\rho_n} \right ) \right )^{n-1}+  \frac{ Cn}{(x \rho_n)^2} F^{n-1}\left (\frac{1}{x\rho_n} \right )  .
		\end{align*}
		Let $a<\frac{1}{2}$ and $n$ be sufficiently large so that $$ F\left (\frac{1}{x\rho_n} \right )< a. $$
		Then $$ \frac{  n}{  \sigma_n ^2} \left (2F\left (\frac{1}{x\rho_n} \right ) \right )^{n-1}<  \frac{  n}{  \sigma_n ^2} (2a)^{n-1}\to 0.$$
		The same argument proves that    $$  \frac{ n}{  \rho_n ^2}   F  ^{n-1}\left (\frac{1}{x\rho_n} \right  )   \leq \frac{  n}{  \rho_n ^2} \left  (2F \left (\frac{1}{x\rho_n}  \right )  \right )^{n-1}\to 0.$$
		For the second and third summand we have that
		\begin{align*}
			&   \alpha_n | b_n-\beta_n |\leq \frac{S_{y\sigma_n} n}{(y\sigma_n)^2} F^{n}\left (\frac{1}{x \rho_n} \right )\left  ( 1+F\left (\frac{1}{{y\sigma_n}}\right )\right )^{n-1}\leq   \frac{C n}{(y\sigma_n)^2} \left (2 F\left (\frac{1}{x \rho_n} \right )\right )^{n-1}  
		\end{align*}
		and
		\begin{align*}
			& b_n|a_n - \alpha_n|\leq |a_n - \alpha_n|\leq  \frac{S_{x \rho_n} n}{(x \rho_n)^2}\left \{F \left (\frac{1}{{x \rho_n}}\right  )\right \}^{n-1}\leq  \frac{Cn}{(x \rho_n)^2}\left \{F \left (\frac{1}{{x \rho_n}}\right  )\right \}^{n-1}.
		\end{align*} 
		Note that convergence to 0 for the latter two terms follows by the same argument as for the first summand.
	\end{proof}
	
	\begin{remark} 
		Assumption \eqref{10} holds if there is $\beta$ such that $F(x) - F(y) \leq \beta (x-y)$ for $x\geq y$   in a neighborhood $(0,\varepsilon)$ of the origin. In fact, if this is the case, we have  $$\sup_{a\geq  \frac{1}{\epsilon} }S_a\leq    \sup_{ a\geq  \frac{1}{\epsilon  }} \beta =\beta ;$$ 
		hence we can take $A=\frac{1}{ \epsilon }$. 
		
	\end{remark} 
	
	\section{Extremal types theorem for generalized Oppenheim expansions}
	In this section we study the behaviour of the sequence  $( Z_n)_{n \geq 1}$   in depth. Our objective is to establish an Extremal Types Theorem akin to the independent case, as outlined in Theorem \ref{extr}. To accomplish this, we must first derive a set of preliminary results. 
	
	\begin{lemma}\label{2} \sl     Assume that condition \eqref{condF} is satisfied. Then for every $k$,  for  every finite increasing sequence of integers $i_1 , \dots,  i_k$  and  for  every finite sequence of numbers $x_{ 1}, \dots, x_{ k} \geq  1$ we have
		$$\left |P(R_{i_ 1}> x_{ 1}, \dots, R_{i_k}> x_{ k})- P(R_{i_1}> x_{ 1}) \cdots P(R_{i_k}> x_{ k})\right |\leq \beta \left (\sum_{j=1}^k \frac{1}{ x^2_{ j}}\prod_{r\not =j} F\left (\frac{1}{x_{ r}}\right )\right ).$$
	\end{lemma}
	\begin{proof} Lemma \ref{1} gives
		\begin{align*}&
			\left |P(R_{i_1}> x_{ 1}, \dots, R_{i_k}> x_{ k})- P(R_{i_1}> x_{ 1}) \cdots P(R_{i_k}> x_{ k})\right |\\&\leq F\left (\frac{1}{ x_{ 1}}\right ) \cdots F\left (\frac{1}{ x_{ k}}\right ) - F\left (\frac{1}{ x_{ 1}+1}\right ) \cdots F\left (\frac{1}{ x_{ k}+1}\right )  .
		\end{align*}
		Now, for any pair of sequences $(a_j)_{1\leq j \leq k}$ and $(b_j)_{1\leq j \leq k}$ we can write
		\begin{equation}\label{13}
			\prod_{j=1}^k a_j - \prod_{j=1}^k b_j = \sum_{r=1}^k\prod_{j=1}^{r-1} a_j(a_r-b_r)\prod_{j=r+1}^{k} b_j
		\end{equation}
		(proof by induction. Recall the convention $\prod_{j\in \emptyset}=1$). This gives
		\begin{align*}&
			F\left (\frac{1}{ x_{ 1}}\right ) \cdots F\left (\frac{1}{ x_{ k}}\right ) - F\left (\frac{1}{ x_{ 1}+1}\right ) \cdots F\left (\frac{1}{ x_{ k}+1}\right ) \\&= \sum_{r=1}^k\prod_{j=1}^{r-1}  F\left (\frac{1}{ x_{ j}}\right ) \left (  F\left (\frac{1}{ x_{ r}}\right )-  F\left (\frac{1}{ x_{ r}+1}\right )\right )\prod_{j=r+1}^{k}   F\left (\frac{1}{ x_{ j}+1}\right ) \leq \beta \left (\sum_{r=1}^k \frac{1}{ x^2_{ r}}\prod_{j\not =r} F\left (\frac{1}{x_{ j}}\right )\right )
		\end{align*}
		since
		$$F\left (\frac{1}{ x_{r}}\right )-  F\left (\frac{1}{ x_{r}+1}\right )\leq   \frac{\beta}{x_{r}(x_{r}+1)}\leq \frac{\beta}{x^2_{r}}.$$
	\end{proof}
	
	\noindent
	The result that follows generalizes Lemma \ref{2} and it can be obtained by applying similar arguments.
	
	\begin{lemma}\label{9} \sl Assume that condition \eqref{condF} is satisfied. Then for every $k$,  every family   $I_1, \dots, I_k$ of disjoint intervals of integers and every $x \geq 1$ we have
		\begin{align*}&
			\left |P\left  (\bigcap_{i \in I_1} \{R_{i }> x \}, \dots, \bigcap_{i \in I_k} \{R_{i }> x \} \right )- P\left  (\bigcap_{i \in I_1} \{R_{i }> x \}\right ) \cdots  P\left  (\bigcap_{i \in I_k} \{R_{i }> x \}\right )\right | \\&\leq   \frac{\beta q}{ x^2 }  F^{q-1}\left (\frac{1}{x }\right ),
		\end{align*}	
		where $q = \sum_{j=1}^k \#I_j. $	 
	\end{lemma}
	
	\noindent
	\noindent

		


	\begin{remark}
		Putting $y = \frac{1}{x}$, the inequality of  Lemma \ref{9} can be also rephrased as
		
		\begin{align}& \nonumber
			\left|P\left (\max_{i \in I_1} \frac{1}{ R_{i }}
			\leq y , \dots,  \max_{i \in I_k}\frac{1}{ R_{i }}
			\leq y \right)- P\left (\max_{i \in I_1}  \frac{1}{ R_{i }}
			\leq y  \right)\cdots  P\left (\max_{i \in I_k}  \frac{1}{ R_{i }}
			\leq y \right)  \right| \\&\label{8}\leq    \beta q y^2      F^{q-1} (y), 
		\end{align}	
		which can be seen as the analogue to Lemma 3.2.2 of \cite{LLR}.
	\end{remark}

	\medskip

	\medskip
	
	\noindent We describe here the construction of some intervals as presented on p. 55 of \cite{LLR}. Precisely, let $k$ be a fixed integer, and for any positive integer $n$, write $n^{\prime}=[n / k]$ (the integer part of $n / k)$. Thus we have $n^{\prime} k \leq n<\left(n^{\prime}+1\right) k$. Divide the first $n^{\prime} k$ integers into $2 k$ consecutive intervals, as follows. For large $n$, let $m$ be an integer, $k<m<n^{\prime}$, and write
	$$
	I_1=\left\{1,2, \ldots, n^{\prime}-m\right\}, I_1^*=\left\{n^{\prime}-m+1, \ldots, n^{\prime}\right\},
	$$
	$$
	I_2=\left\{n^{\prime}+1, \ldots, 2n^{\prime}-m\right\}, I_2^*=\left\{2n^{\prime}-m+1, \ldots, 2n^{\prime}\right\},
	$$
	$I_3, I_3^*, \ldots, I_k, I_k^*$ being defined similarly, alternatively having length $n^{\prime}-m$ and $m$. Finally, write
	$$
	I_{k+1}=\left\{(k-1) n^{\prime}+m+1, \ldots, k n^{\prime}\right\}, \quad I_{k+1}^*=\left\{k n^{\prime}+1, \ldots, k n^{\prime}+m\right\} .
	$$
	Observe that $I_{k+1}, I_{k+1}^*$ are defined differently from $I_i, I_j^*$ for $j \leq k$.
	
	\medskip
	
	\noindent The notation introduced above is used in the result that follows which is motivated by Lemma 3.3.1 of \cite{LLR}.

	\begin{proposition}\label{prop3}\sl  Assume that condition \eqref{condF} is satisfied and let $x\leq 1$.  Then  the following inequalities hold true. 
		\begin{enumerate}
			\item[(i)] $$0 \leq P\left(\displaystyle \bigcap_{j=1}^{k}\left\{\max_{i\in I_j}\dfrac{1}{R_i}<x\right\}\right)- P\left(\max_{1 \leq j\leq n}\dfrac{1}{R_j}<x\right)\leq \sum_{j =1}^{k+1}P\left( \max_{i\in I_j}\dfrac{1}{R_i}<x <\max_{i\in I^*_j}\dfrac{1}{R_i}\right);$$
			
			\item[(ii)]	\begin{align*}&
				\left|P \left(\displaystyle \bigcap_{j=1}^{k}\left\{\max_{i\in I_j}\dfrac{1}{R_i}<x\right\}\right)- \prod_{j=1}^{k }P\left(\max_{i\in I_j}\dfrac{1}{R_i}<x\right)\right|\leq \beta (k-1)(n^\prime - m)x^2F^{ k(n^\prime - m)-1}(x) ;
			\end{align*}
			
			\item[(iii)] \begin{align*}
				&\left| \prod_{j=1}^{k }P\left(\max_{i\in I_j}\dfrac{1}{R_i}<x\right)-P^k \left(\max_{i\in I_1}\dfrac{1}{R_i}<x\right)\right|\leq \sum_{k=1}^k \left|P\left(\max_{i\in I_j}\dfrac{1}{R_i}<x\right)-P  \left(\max_{i\in I_1}\dfrac{1}{R_i}<x\right)\right|\\
				& {\leq \beta k (n^\prime - m)x^2 F^{(n^\prime - m)-1}(x);}\\
			\end{align*}
			\item[(iv)]\begin{align*}&\left|P^k \left(\max_{i\in I_1}\dfrac{1}{R_i}<x\right)- P^k \left(\max_{1\leq i\leq n^\prime }\dfrac{1}{R_i}<x\right)\right|\leq k P\left( \max_{i\in I_1}\dfrac{1}{R_i}<x <\max_{i\in I^*_1}\dfrac{1}{R_i}\right).\end{align*}
		\end{enumerate}
		
	\end{proposition}
	\begin{proof}
		For (i), observe that 
		$$
		\bigcap_{j=1}^k\left\{\max _{  i\in I_j} \frac{1}{R_ i}<x\right\} \supset\left\{\max _{1 \leq i \leq n} \frac{1}{R_i}<x\right\},
		$$
		and their difference implies that 
		$$\max _{i \in I_j} \frac{1}{R_i}<x<\max _{ i \in I_j^*} \frac{1}{R_i}$$ for some $j \leq k$
		or otherwise 
		$\frac{1}{R_j} \leq x$ for $1<j \leq k n^{\prime}$ but $\frac{1}{R_j}>x$ for some $j=k n^{\prime}+1, \ldots, k\left(n^{\prime}+1\right)$ which implies that
		$$
		\max _{j \in I_{k+1}} \frac{1}{R_j}<x<\max _{j \in I_{k+1}^*} \frac{1}{R_j}  
		$$
		since $m>k \Rightarrow k\left(n^{\prime}+1\right)<k n^{\prime}+m$.   
		
		\medskip
		\noindent
		Part (ii) follows as a direct consequence of \eqref{8}. For part (iii) we start by  employing equation \eqref{13}, to write
		\begin{align*}&
			\prod_{j=1}^{k }P\left(\max_{i\in I_j}\dfrac{1}{R_i}<x\right)-P^k \left(\max_{i\in I_1}\dfrac{1}{R_i}<x\right) \\&
			= \sum_{r=1}^k \prod_{j=1}^{r-1} P\left(\max_{i\in I_j}\dfrac{1}{R_i}<x\right) \left\{P\left(\max_{i\in I_r}\dfrac{1}{R_i}<x\right)-P \left(\max_{i\in I_1}\dfrac{1}{R_i}<x\right)\right\}\cdot\\&\cdot \prod_{j=r+1}^{k}  P \left(\max_{j\in I_1}\dfrac{1}{R_j}<x\right);
		\end{align*}		
		hence, since the product terms are not larger than 1,	
		\begin{align*}&\left| \prod_{j=1}^{k }P\left(\max_{i\in I_j}\dfrac{1}{R_i}<x\right)-P^k \left(\max_{i\in I_1}\dfrac{1}{R_i}<x\right)\right|\leq \sum_{j=1}^k \left|P\left(\max_{i\in I_j}\dfrac{1}{R_i}<x\right)-P  \left(\max_{i\in I_1}\dfrac{1}{R_i}<x\right)\right|\\&=\sum_{j=2}^k \left|P\left(\max_{i\in I_j}\dfrac{1}{R_i}<x\right)-P  \left(\max_{i\in I_1}\dfrac{1}{R_i}<x\right)\right|.
		\end{align*}	
		For every $j=1, \dots, k$ we have, by  Lemma \ref{1} (ii),
		$$F^{n^\prime -m}\left (\frac{x}{x+1}\right )\leq P\left(\max_{i\in I_j}\dfrac{1}{R_i}<x\right)\leq F^{n^\prime -m}\left (x\right );$$
		hence, for every $j=2, \dots, k$
		\begin{align*}&
			\left|P\left(\max_{i\in I_j}\dfrac{1}{R_i}<x\right)-P  \left(\max_{i\in I_1}\dfrac{1}{R_i}<x\right)\right|\leq F^{n^\prime -m}\left (x\right )-F^{n^\prime -m }\left (\frac{x}{x+1}\right )\\
			&\leq \beta(n^\prime -m) \left(x -\frac{x}{x+1} \right )F^{n^\prime -m-1}\left (x\right )\leq \beta(n^\prime -m)  x^2F^{n^\prime -m-1}\left (x\right ).
		\end{align*}
		which leads to the desired result.
		
		\medskip
		\noindent
		Finally for part (iv), by the same argument as in (iii), we have
		\begin{align*}&
			\left|P^k \left(\max_{i\in I_1}\dfrac{1}{R_i}<x\right)- P^k \left(\max_{1\leq i\leq n^\prime }\dfrac{1}{R_i}<x\right)\right|\leq k \left|P  \left(\max_{i\in I_1}\dfrac{1}{R_i}<x\right)- P  \left(\max_{1\leq i\leq n^\prime }\dfrac{1}{R_i}<x\right)\right|,
		\end{align*}		
	\end{proof}	
	
	\begin{proposition}\label{prop4}\sl  Assume that condition \eqref{condF} is satisfied and consider a sequence of real numbers $(u_n)_{n\geq 1}$ such that $u_n\to 0$ as $n\to\infty$. Then,
		$$\lim_{n \to \infty}\left\{P\left(\max_{1 \leq j\leq n}\dfrac{1}{R_j}<u_n\right)- P^k \left(\max_{1\leq j\leq n^\prime }\dfrac{1}{R_j}<u_n\right)\right\}=0$$
	\end{proposition}	
	\begin{proof}
		Summing the right hand members of the inequalities in the above Proposition \ref{prop3} we obtain
		\begin{align*}&
			\left|P\left(\max_{1 \leq j\leq n}\dfrac{1}{R_j}<u_n\right)- P^k \left(\max_{1\leq j\leq n^\prime }\dfrac{1}{R_j}<u_n\right)\right|\\
			&\leq  \sum_{j =1}^{k+1}P\left( \max_{i\in I_j}\dfrac{1}{R_i}<u_n <\max_{i\in I^*_j}\dfrac{1}{R_i}\right)\\
			& +\beta k(n^\prime - m)u_n^2F^{ k(n^\prime - m)-1}(u_n)+\beta (k-1) (n^\prime - m)u_n^2 F^{(n^\prime - m)-1}(u_n) \\
			&+k P\left( \max_{i\in I_1}\dfrac{1}{R_i}<u_n <\max_{i\in I^*_1}\dfrac{1}{R_i}\right).
		\end{align*}
		Obviously the summand
		$$	 \beta k(n^\prime - m)u_n^2F^{ k(n^\prime - m)-1}(u_n)+\beta (k-1) (n^\prime - m)u_n^2 F^{(n^\prime - m)-1}(u_n) $$
		goes to 0 as $n \to \infty.$ So it remains to bound conveniently the summands
		$$P\left( \max_{i\in I_j}\dfrac{1}{R_i}<u_n <\max_{i\in I^*_j}\dfrac{1}{R_i}\right), \qquad j=1 , \dots, k+1.$$
		Note that we can write
		\begin{align*}
			&P\left( \max_{i\in I_j}\dfrac{1}{R_i}<u_n <\max_{i\in I^*_j}\dfrac{1}{R_i}\right) =  P\left( \left \{\max_{i\in I_j}\dfrac{1}{R_i}<u_n\right \}  \cap \left \{\max_{i\in I^*_j}\dfrac{1}{R_i}>u_n\right \}\right)\\
			& = P\left ( \max_{i\in I_j}\dfrac{1}{R_i}<u_n\right ) - P\left(\max_{i\in I_j}\dfrac{1}{R_i}<u_n ,  \, \max_{i\in I^*_j}\dfrac{1}{R_i}\leq u_n\right)\\
			& =  P\left ( \bigcap_{i\in I_j}\left \{\dfrac{1}{R_i}<u_n\right \}\right ) - P\left(\bigcap_{i\in I_j}\left \{\dfrac{1}{R_i}<u_n\right \} ,  \,\bigcap_{i\in I^*_j}\left \{\dfrac{1}{R_i}\leq u_n\right \}\right)\\
			&\leq F^{n^\prime -m}\left (u_n\right )-F^{n^\prime -m}\left (\dfrac{u_n}{u_n+1}\right )F^{m}\left (\dfrac{u_n}{u_n+1}\right )
		\end{align*}
		Thus,

		\begin{align*}&
			\sum_{j =1}^{k+1}P\left( \max_{i\in I_j}\dfrac{1}{R_i}<u_n <\max_{i\in I^*_j}\dfrac{1}{R_i}\right)+k P\left( \max_{i\in I_1}\dfrac{1}{R_i}<u_n <\max_{i\in I^*_1}\dfrac{1}{R_i}\right)\\ 
			&\leq (2k+1)\left(F^{n^\prime -m}\left (u_n\right )-F^{n^\prime -m}\left (\dfrac{u_n}{u_n+1}\right )F^{m}\left (\dfrac{u_n}{u_n+1}\right )\right).
		\end{align*}
		which tends to 0 as $n\to \infty$.   
	\end{proof}

	\begin{theorem}\sl  
		\label{extr}{Assume that condition \eqref{condF} is satisfied and} that there exists constants $a_n\to \infty$ and $b _n  $ such that
		$$P\left (a_n \left ( \max_{1 \leq i \leq n}\dfrac{1}{R_i}-b _n  \right )\leq x\right )$$
		converges to a non degenerate distribution function $G$.  Then $G$ has one of the threee extreme values form.
	\end{theorem}
	\begin{proof}
		We want to prove that $G$ is an extreme value distribution. According to Theorem 1.3.1 in \cite{LLR} (p.8), this follows if we prove that for each $k=1,2,3,\ldots$
		\[
		P\left(\max_{1\leq j\leq n} \dfrac{1}{R_j} \leq \dfrac{x}{a_{nk}}+b_{nk}\right) \xrightarrow{\text{w}}G^{1/k}(x)
		\]
		where $\xrightarrow{\text{w}}$ denotes convergence at continuity points of the limiting function. By Proposition \ref{prop4}, the latter expression is equivalent to 
		\[
		P\left( \max_{1\leq j\leq nk}\dfrac{1}{R_j}\leq\dfrac{x}{a_{nk}}+b_{nk}\right) - P^k\left( \max_{1\leq j\leq n}\dfrac{1}{R_j}\leq\dfrac{x}{a_{nk}}+b_{nk}\right) \to 0\quad \mbox{as} \quad n\to\infty
		\]
		for all integers $k$. Observe that $k=1$ is true by assumption, thus we need to prove the relation for $k=2,3, \ldots.$ This follows immediately by applying Proposition \ref{prop4} for $nk$ instead of $n$.
	\end{proof}
	
	\section{Discussion}
	Lemma \ref{9}, Proposition \ref{prop3}, Proposition \ref{prop4} and Theorem \ref{extr} were motivated by the work presented in Chapter 3 of \cite{LLR} for stationary sequences of random variables. For obtaining the corresponding results for stationary sequences, the authors assumed the validity of the so-called Leadbetter condition: Let $\left(u_n\right)_{n\geq 1}$ be a sequence of real numbers and consider integers such that $$1 \leq i_1<i_2<\cdots<i_p<j_1 \cdots<j_q \leq n, j_1-i_p \geq \ell.$$ A sequence of random variables $(X_n)_{n\geq 1}$ is said to satisfy condition $\Delta(u_n)$ if 
	$$
	\begin{aligned}
		& \mid P\left(X_{i_1} > u_n, \cdots, X_{i_p} > u_n, X_{j_1} > u_n, \cdots, X_{j_q} > u_n\right) \\
		& \quad-P\left(X_{i_1} > u_n, \ldots, X_{i_p} > u_n\right) P\left(X_{j_1} > u_n, \ldots, X_{j_q} > u_n\right) \mid \leq \alpha_{n, \ell},
	\end{aligned}
	$$
	where $\alpha_{n, \ell}$ is nonincreasing in $\ell$ and such that
	$$
	\lim _{n \rightarrow \infty} \lim _{n \rightarrow \infty} \alpha_{n, \ell_n}=0
	$$
	for some sequence $\left(\ell_n\right)_n \rightarrow \infty$ with $\frac{\ell_n}{n} \rightarrow 0$. 
	
	\medskip
	
	\noindent It can be proven that the random variables $(R_n)_{n\geq 1}$ satisfy a similar condition as it can be seen in the Theorem that follows. It is worth noting that this result served as a motivation for studying the particular topic.
	
	\begin{theorem}\sl
		If \eqref{condF} is satisfied, then condition $\Delta\left(u_n\right)$ holds for the random variables $\left(R_n\right)_{n\geq 1}$ for every sequence $\left(u_n\right)_{n\geq 1}$ such that $\displaystyle\lim _{n \rightarrow \infty} u_n=\infty$.
	\end{theorem}
	\begin{proof}
		For brevity denote
		$$
		\begin{gathered}
			P\left(R_{i_1}>u_n, \ldots, R_{i_p}>u_n, R_{j_1}>u_n, \ldots, R_{j_q}>u_n\right)=\rho_{p+q}(n) ; \\
			P\left(R_{i_1}>u_n, \ldots, R_{i_p}>u_n\right)=\rho_p(n) ; \quad P\left(R_{j_1}>u_n, \ldots, R_{j_q}>u_n\right)=\rho_q(n) .
		\end{gathered}
		$$
		Then, from Theorem \ref{6}  we have 
		$$
		\begin{aligned}
			& \left|\rho_{p+q}(n)-\rho_p(n) \rho_q(n)\right|= \left |\left(\rho_{p+q}(n)-F^{p+q}\left(\frac{1}{u_n}\right)\right)-\left (\rho_p(n)-F^p\left(\frac{1}{u_n}\right)\right )\left ( \rho_q(n)-F^q\left(\frac{1}{u_n}\right)\right )  \right.\\
			& \left.-F^p\left(\frac{1}{u_n}\right)\left(\rho_q(n)-F^q\left(\frac{1}{u_n}\right)\right)-F^q\left(\frac{1}{u_n}\right)\left (\rho_p(n)-F^p\left(\frac{1}{u_n}\right)\right)  \right | \\
			& \leq \frac{\beta(p+q)}{u_n^2} F^{p+q-1}\left(\frac{1}{u_n}\right)+ \frac{\beta^2 pq}{u_n^4}F^{p+q-2}\left(\frac{1}{u_n}\right)
			+\frac{\beta p}{u_n^2} F^{p+q-1}\left(\frac{1}{u_n}\right)+\frac{\beta q}{u_n^2} F^{p+q-1}\left(\frac{1}{u_n}\right)   \\
			\\& \leq \frac{2\beta(p+q)}{u_n^2}  +  \frac{\beta^2 pq}{u_n^4} 
			\rightarrow 0.
		\end{aligned}
		$$
	\end{proof}
	
	\begin{remark}
		The result ensures that long-range dependence is sufficiently weak so as not to affect the asymptotics on an exteme value analysis.
	\end{remark}
	
	\noindent The next remark concerns Theorem \ref{extr}. 
	
	\medskip
	\begin{remark}  
		Regrettably, we currently lack an example demonstrating the application of {\color{red} Theorem \ref{extr}}. This observation bolsters the conjecture that the limit of every normalized sequence  $a_n  \left ( \max_{1 \leq i \leq n}\dfrac{1}{R_i}-b_n\right ) $  (with $a_n \to \infty$) either doesn't exist or is degenerate (refer to Corollary \ref{11} and Remark \ref{15}).  We intend to delve deeper into this issue in the future.
	\end{remark}

\end{document}